\documentclass[a4paper,leqno,11pt,final]{amsproc}
\usepackage{amssymb,latexsym,amsxtra,amscd,amsfonts,amsthm,amsmath,verbatim}

\newcommand\QQ{\mathbb{Q}}

\newcommand\codim{\mathrm{codim}}
\newcommand\Hom{\mathrm{Hom}}
\newcommand\Spec{\mathrm{Spec}}

\numberwithin{equation}{section}

\newtheorem{theorem}[equation]{Theorem}
\newtheorem{example}[equation]{Example}

\begin{document}


\title[Methods for Computing Normalisations of Affine Rings]
{Methods for Computing Normalisations\\ of Affine Rings}

\author{\sc A. Taylor}

\thanks{The research of the author was paritally supported the National Security Agency.}

\address{Department of Mathematics,
Rutgers University, Piscataway, NJ 08854}

\email{ataylor@math.rutgers.edu}

\begin{abstract}
Our main purpose is to give multiple examples for 
using the available implementations for computing the 
normalization of an affine ring, computing the minimial 
generators of the normalization as an algebra over the 
original ring and integral closures of ideals.  Some 
such examples have been published for {\sc Singular}, but 
not for {\sc Macaulay 2} and we present both in this paper.  
We also briefly describe the implementations.  
\end{abstract}

\maketitle


\section{Introduction}
In this note we describe the available methods for computing the 
normalization of affine rings, that is reduced, finite, $k$-algebras, 
where $k$ is a computable field, and the implications for computing the 
integral closure of ideals.  Our goal is to keep it simple.  
Much more is known for computing the normalization with assumptions on 
the ring or the ideal, as well as for the normalization of discrete 
valuation rings.  We only address what can be done in the most 
general setting of affine rings.  

Our main goal is to give examples so that users can fully exploit the 
implementations.  We also briefly describe which algorithms have been 
implemented, some implementation issues and how to use these implementations 
to compute the integral closure of small ideals.  All 
computations listed in this paper were done on a Pentium III, 600 MHZ 
machine with 256 MB of RAM.  

At the time of this writing there are two implementations of normalization.  
One implementation is in {\sc Macaulay 2}~\cite{GS} and  the other is in 
{\sc Singular}~\cite{GPS01}.  Both implementations are based 
on the information given in de Jong's paper~\cite{DJ}.  The philosophical 
approach of which is seen in Vasconcelos' papers~\cite{V2} and \cite{V4} 
and a more detailed analysis of this approach is given in~\cite{DGJP}. 
A brief description as well as other examples for {\sc Singular} are given 
in~\cite{GM1}. Finally, there will be chapter on normalization, computing 
it in {\sc Singular}, and related issues in the upcoming book of Greuel and 
Pfister~\cite{GP1}.  All of what is published is based on the implementation 
in {\sc Singular}, where as this is the first document detailing what is 
done in {\sc Macaulay 2} as well.

Let $A = R/I$ be an affine ring.  In the following we use the language of 
Matsumura~\cite{M}.  If $A$ is a domain and $Q$ its ring of fractions then 
the \emph{integral closure} $\overline{A}$ of $A$ is the set of all elements 
in $Q$ that are integral over $A$ and $A$ is \emph{integrally closed} if 
$A = \overline{A}$.  An affine ring $A$ is \emph{normal} if $A_p$ is an 
integrally 
closed domain at every prime $p\in \Spec(A)$. If $A$ is normal it is 
well known that $A\cong A/P_1\times \cdots A/P_t$ where $P_1,\ldots P_t$ are the 
minimal primes for $A$ and $A/P_i$ is an integrally closed domain for 
$1\leq i\leq t$.  Since the implementations deal with affine rings in general 
this decomposition is often referred to as the \emph{normalization} of $A$ 
rather than the integral closure.  {\sc Macaulay 2} refers to the computation 
as integral closure where as {\sc Singular} and all of the publications 
associated to its implementation call the computation normalization.  To 
emphasize that we do not require the ring to be a domain, we also refer 
to the computation as normalization.  

Serre's criterion for normality 
implies that an affine ring $A$ is normal if and 
only if it is both $S_2$ and $R_1$~\cite[Theorem 23.8]{M}. An affine ring 
is $S_2$ if $A_p$ is Cohen-Macaulay for $p\in \Spec(A)$ and  
$\codim(p)\leq 2$, and has depth greater than or equal to 2 for 
$p\in \Spec(A)$ and $\codim(p)> 2$.  An affine ring is $R_1$ if 
it is regular in codimension one.  

The philosophy of the approach is 
to enlarge $A$ recursively inside of 
$Q$ until the normalization of $A$ is obtained.  The extension Vasconcelos 
uses in~\cite{V2}, \cite{V4} is $\Hom(\Hom(L,A),\Hom(L,A))$, where 
$L$ is the Jacobian ideal of $A$.   
The following key theorem in de Jong's paper~\cite{DJ} is originally due to 
Grauert and Remmert~\cite[pp. 220-221]{GR71}, \cite[pp. 125-127]{GR84} and 
describes the ring de Jong uses for this 
extension.  The \emph{non-normal locus} of $A$ is the set $NNL\subseteq \Spec(A)$ 
such that for $p\in NNL$, $A_p$ is not an integrally closed domain.

\begin{theorem}\label{dJTest} Let $A = R/I$ be an affine ring and $J$ an 
ideal of $A$.  Assume that the ideal $J$ contains a non-zero divisor, 
and has the following property:
$$ NNL \subseteq V(J), $$
where $V(J) = \{p\in \Spec(A) | J\subset p\}$ denotes, as usual, the 
zero set of $J$.  Suppose moreover that $J$ has the property
\begin{equation}\label{eq:dJTest} \Hom_A(J,J) = \Hom_A(J,A)\cap\overline{A}.
\end{equation}
Then one has the following normality criterion:
$$ A = \Hom_A(J,J) \hbox{ if and only if } A \hbox{ is normal. }$$
\end{theorem}

The implementations in both {\sc Macaulay 2} and {\sc Singular} are split 
into two phases.  First they find an ideal $J$ such that the non-normal 
locus of $A$ is contained in $V(J)$ and satisfying 
Equations~(\ref{eq:dJTest}).  In the second phase $\Hom(J,J)$ is computed 
as a ring.  

For the first phase, Theorem 2.2 in de Jong's paper~\cite{DJ} establishes 
that any radical ideal containing a non-zero divisor satisfies 
Equation (\ref{eq:dJTest}).  Also, it is well known 
that if the Jacobian ideal $L$ of $A$ is not contained in $Q\in\Spec(A)$ then 
$A_Q$ is regular, so the non-normal locus of $A$ is contained in 
$V(L) = V(\sqrt{L})$.  Hence $\sqrt{L}$ satisfies Theorem~\ref{dJTest}.  
Computing $L$ can be very complex in terms of space and time, so 
sometimes the programs may use $\sqrt{f}$ for $f\in L$ non-zero which 
also satisfies the hypotheses of Theorem~\ref{dJTest}.  

 For the second phase, there are two facts that are very helpful.  The first, 
is the well known fact that if $f$ is a non-zero divisor in $A$,
then $(fJ:J) = f\Hom(J,J)$.  The second, is the presentation of 
$\Hom(J,J)$ as a ring that is due to Catanese~\cite{CA} 
and is given in de Jong's paper~\cite[Theorem 3.1]{DJ}.  Let 
$v_1,\ldots,v_n$ be a module basis for $\Hom(J,J)$.  Map 
$S = k[X_1,\ldots,X_n]$ onto $\Hom(J,J)$ by sending $X_i$ to $v_i$ for 
$1\leq i\leq n$.  Then as a ring $\Hom(J,J)$ is defined by two sets of 
relations.  There are linear equations which are 
the syzygies for the module, $\sum \alpha_{i}X_i = 0$ and quadratic 
equations which come from the 
fact that $\Hom(J,J)$ is a ring so there exists $\beta_{ijk}$, $1\leq k\leq n$ 
such that $v_iv_j = \sum\beta_{ijk}v_k$.  It is a theorem of Catanese 
stated clearly in \cite{DGJP} and \cite{DJ} that the ideal generated by 
these two sets of relations is the full kernel of the map, thus inducing an 
isomorphism.

The implementations in {\sc Macaulay 2} and {\sc Singular} were completed 
independently, but the pseudo-code for both is very close to the 
presentation in~\cite{DGJP} so we will not include it here.  

The casual user looking at the code in either program will see 
relatively complicated programs which is due to 
the recursive nature of the algorithm.  The program must track the proper 
extensions and in the case of reduced rings, it must track the splitting
of the ring into domains.  Often for an ideal $J$ as in the theorem, while 
$\Hom(J,J)$ is a proper extension, it is not the normalization yet so 
the process must be repeated many times.  It is unknown, given the ring, 
the number of extensions that must be computed.  
Moreover, experimental experience suggests that there are choices available 
that could reduce the number of extensions that need to be computed, but knowing 
when to use them is not always a priori clear.  There is a new algorithm of 
Vasconcelos' for which the number of recursions is bounded, but this algorithm 
is not yet implemented in any of the systems~\cite{V3}.

\section{Examples}
For the examples we will focus on a hypersurface given to us by Craig 
Huneke, a union of straight lines given in the {\sc Singular} example files and 
an example of Huneke's given in each of~\cite{DGJP}, \cite{V1}, \cite{V2}, 
and \cite {V4}.  Also the first and third examples were chosen because one 
satisfies Serre's condition $S_2$ but does not satisfy Serre's condition 
$R_1$ and the other is the opposite.  

For each of the examples we give the input and output of 
both programs so that the user can simply type in the same information and  
see the same results.  Occasionally, due to the width of the text we have 
had to slightly alter the output for it to fit in the space allotted. 
We hope that this is also enough information for the 
user to then perform similar computations on their own examples.  Also, when 
possible we include the computation times in CPU seconds.  If no time is 
listed it is because the time is negligible.  We include the times because they 
are indicative of what can happen on larger examples and because they give more 
meaning to the times listed after Example~\ref{complicated}.  However, when 
taken by themselves the times may not have much meaning since even on the same 
machine times can vary up to a few seconds.  

\begin{example}\label{ex1:hypersurface} Let 
$A = \QQ[x,y,z]/(x^6-z^6-y^2z^4)$. This ring is not $R1$, but it is 
$S2$. Therefore, computing the normalization in the example forms the 
desingularization of $A$ in codimension one. The computation 
in {\sc Macaulay 2} follows.
\begin{verbatim}
i1 : A = QQ[x,y,z]/ideal(x^6-z^6-y^2*z^4)

o1 = A

o1 : QuotientRing

i2 : time integralClosure (A, Variable => symbol a)
     -- used 1.67 seconds

                  
                 QQ [a , a , x, y, z]
                      7   6
o2 = ---------------------------------------------
       2                    2         2    2    2
     (x  - a z, a x - a z, a  - a x, a  - y  - z )
            6    6     7    6    7    7

o2 : QuotientRing
\end{verbatim}
Computing the normalization of an affine ring requires new 
variables unless the ring is already normal.  {\sc Macaulay 2} gives the user 
the opportunity to specify which letter to use for this new variable.  
This is done via
\begin{verbatim}
Variable => symbol a.
\end{verbatim}
If this is not included, then {\sc Macaulay 2} uses $w$ by default.  
Using a letter to create the normalization which is a variable in the definition 
of $R$ can lead to problems with using the 
normalization, so we recommend supplying a variable not used to define $R$.  

For {\sc Singular}, we must first load the normalization library. 
{\sc Singular} then lists the libraries loaded and we use elipses 
to indicate that more is printed on that line, but we have removed 
it to keep from having to make multiple lines due to the size of the 
paper. Also, 
we use nor to name the output of normal(I) to match the comments that 
singular outputs as part of this function.  We then follow the other 
commands they list for continuity.

\begin{verbatim}
> LIB "normal.lib";
// ** loaded /usr/local/Singular/2-0-3/LIB/normal.lib ...
// ** loaded /usr/local/Singular/2-0-3/LIB/hnoether.lib ... 
// ** loaded /usr/local/Singular/2-0-3/LIB/primitiv.lib ... 
// ** loaded /usr/local/Singular/2-0-3/LIB/random.lib ... 
// ** loaded /usr/local/Singular/2-0-3/LIB/matrix.lib ... 
// ** loaded /usr/local/Singular/2-0-3/LIB/ring.lib ... 
// ** loaded /usr/local/Singular/2-0-3/LIB/inout.lib ... 
// ** loaded /usr/local/Singular/2-0-3/LIB/presolve.lib ... 
// ** loaded /usr/local/Singular/2-0-3/LIB/elim.lib ... 
// ** loaded /usr/local/Singular/2-0-3/LIB/primdec.lib ... 
// ** loaded /usr/local/Singular/2-0-3/LIB/sing.lib ... 
// ** loaded /usr/local/Singular/2-0-3/LIB/poly.lib ... 
// ** loaded /usr/local/Singular/2-0-3/LIB/general.lib ... 
> ring R = 0,(x,y,z),dp;
> ideal I= x6-z6-y2z4;
> list nor = normal(I);

// 'normal' created a list of 1 ring(s).
// nor[1+1] is the delta-invariant in case of choose=wd.
// To see the rings, type (if the name of your list is nor):
     show( nor);
// To access the 1-st ring and map (similar for the others), 
     type: def R = nor[1]; setring R;  norid; normap;
// R/norid is the 1-st ring of the normalization and
// normap the map from the original basering to R/norid
> timer=1;
//used time: 0.82 sec
> show(nor);

// list, 1 element(s):
[1]:
   // ring: (0),(T(1),T(2),T(3),T(4),T(5)),
                           (a(1,1,1,1,1),dp(5),C);
   // minpoly=0
// objects belonging to this ring:
// normap               [0]  ideal, 3 generator(s)
// norid                [0]  ideal, 4 generator(s)
> setring S; norid;
norid[1]=T(4)^2-T(1)*T(5)
norid[2]=T(1)*T(4)-T(3)*T(5)
norid[3]=T(2)^2+T(3)^2-T(5)^2
norid[4]=T(1)^2-T(3)*T(4)
\end{verbatim}

The the output of the command show(nor); states that the normalization 
of $R/I$ requires only one ring because $R/I$ was a domain.  If we 
denote the normalization of $R/I$ by $S/L$ then the next line tells us that 
$S = \QQ[T(1),T(2),T(3),T(4),T(5)]$, each variable has degree 1 and the 
ordering is degree reverse lexicographic.  The next two lines tell us that 
in {\sc Singular} there are two other objects assigned to this ring.  The 
first defines the natural map from $R/I$ into $S/L$ and the second is the 
the defining ideal $L$ of the integral closure.  To see generators of $L$ 
we execute the two commands setring S; and norid;.  
\end{example}

The information stored in the outputs of the two programs 
is fundamentally the same, but the presentation of the 
output in the two programs is a little 
different.  {\sc Macaulay 2} preserves the fact that
$A\subset \Hom(J,J)$ and thus computes a presentation that includes 
those variables from $A$ which still contribute to a minimal 
presentation.  In contrast, {\sc Singular} uses all new variables and in 
this example $T(1) = x, T(2) = y, T(3) = z, T(4) = a_6, T(5) = a_7$.  The 
fractions computed in 
Example~\ref{ex:fractions} can be used to find this correspondence.  
{\sc Macaulay 2} allows the user to choose the name of the 
new variables, but {\sc Singular} does not have such an option.  
When {\sc Singular} computes $normal(I)$ it gives 
information about the number of rings in the decomposition and stores 
this information in the list we called nor above and then a few 
commands are used to display this information.  In this example 
there is one ring, but in the next example 
when the input is not a domain we are told there are 3 rings, so we do the 
display information three times.  In {\sc Macaulay 2} for 
a domain, the domain that is the normalization is displayed 
unless we suppress it using a semicolon.  When the 
ring is not a domain, a sequence of normal 
rings, such that the direct sum of these rings is the 
normalization of the reduced ring is displayed. 

\begin{example}\label{ex2:notAdomain}  Let $A$ be a union of lines, 
$A = \QQ[x,y,z]/((x-y)(x-z)(y-z))$.
\begin{verbatim}
i1 : A = QQ[x,y,z]/ideal((x-y)*(x-z)*(y-z))

i2 : time integralClosure(A,Variable => V)
     -- used 0.47 seconds

      QQ [V , y, z]
           0           R
o2 = (-------------, -----)
          2          y - z
         V  + V
          0    0

o2 : Sequence
\end{verbatim}

Before proceeding to the {\sc Singular} example we consider the 
presentation of the first ring.  The output of {\sc Macaulay 2} will always give 
a normal ring, but not necessarily a domain.  Here it gives 
$$\frac{\QQ[V_0,y,z]}{V_0^2+V_0}\cong\frac{\QQ[V_0,y,z]}{V_0+1}\oplus\frac{\QQ[V_0,y,z]}{V_0}$$
which is normal.  In contrast, {\sc Macaulay 2} computes the normalization of 
$\QQ[V_0,y,z]/(V_0^2+V_0)$ as the isomorphic presentation.  
 
Using {\sc Singular}, the output of $normal(I)$
reveals that there are three rings. What follows reveals that each 
of the three rings is isomorphic to $\QQ[x,y]$.  This is the information we 
gained from {\sc Macaulay 2} in an isomorphic format.  Since we ran this 
example in the same session as the previous example, we do not need 
\begin{verbatim}
> LIB "normal.lib";
\end{verbatim}
this time.  However, if the {\sc Singular} session is restarted then this 
command is needed to get the same results.

\begin{verbatim}
> ring R=0,(x,y,z),dp; ideal I=(x-y)*(x-z)*(y-z); 
> list nor=normal(I);

// 'normal' created a list of 3 ring(s).
// nor[3+1] is the delta-invariant in case of choose=wd.
// To see the rings, type (if the name of your list is nor):
     show( nor);
// To access the 1-st ring and map (similar for the others), 
     type: def R = nor[1]; setring R;  norid; normap;
// R/norid is the 1-st ring of the normalization and
// normap the map from the original basering to R/norid
> timer=1;
> show(nor);
// list, 3 element(s):
[1]:
   //ring: (0),(T(1),T(2)),(dp(2),C);
   // minpoly = 0
// objects belonging to this ring:
// normap               [0]  ideal, 3 generator(s) 
// norid                [0]  ideal, 1 generator(s)
[1]:
   //ring: (0),(T(1),T(2)),(dp(2),C);
   // minpoly = 0
// objects belonging to this ring:
// normap               [0]  ideal, 3 generator(s) 
// norid                [0]  ideal, 1 generator(s)
[1]:
   //ring: (0),(T(1),T(2)),(dp(2),C);
   // minpoly = 0
// objects belonging to this ring:
// normap               [0]  ideal, 3 generator(s) 
// norid                [0]  ideal, 1 generator(s)
> def S1=nor[1]; def S2=nor[2]; def S3=nor[3];
> setring S1; norid; setring S2; norid; setring S3; norid;
norid[1]=0
norid[1]=0
norid[1]=0
\end{verbatim}
Since there are three rings whose direct sum forms the normalization 
of $A = R/I$ the output of show(nor); gives the defining rings for all 
three and we see they are the same.  Then we ask for the defining ideal 
and get that it is zero, so 
$$\overline{A} \cong \QQ[T(1),T(2)]\oplus\QQ[T(1),T(2)]\oplus\QQ[T(1),T(2)]$$ 
which is essentially the same information we got from {\sc Macaulay 2}.
\end{example}

\begin{example}\label{ex3:Wolmer}
This example was chosen to because this 
ring has an isolated singularity, so it is $R_1$, but is not $S_2$ in 
contrast to the first example.  Thus the normalization of this ring 
is its $S_2$-ification.  This example also takes a little longer than the 
previous examples. 
Let $I$ be the radical of the ideal generated by 
\begin{equation*}
\begin{split}
ab^3c+bc^3d+a^3be+cd^3e+ade^3, \\
a^2bc^2+b^2cd^2+a^2d^2e+ab^2e^2+c^2de^2,\\
a^5+b^5+c^5+d^5-5abcde+e^5
\end{split}
\end{equation*}
We use a few new commands in {\sc Macaulay 2} this 
time to make the output more presentable for this paper, however these are 
not needed when a user runs {\sc Macaulay 2} in emacs since the output 
is reasonable with wrap-around turned off.  Even with these new commands we 
altered the presentation of the output a little for both programs to make 
them more readable in this paper.
\begin{verbatim}
i1 : R = QQ[a..e];

i2 : time I = radical(ideal(
a*b^3*c+b*c^3*d+a^3*b*e+c*d^3*e+a*d*e^3,
a^2*b*c^2+b^2*c*d^2+a^2*d^2*e+a*b^2*e^2+c^2*d*e^2,
a^5+b^5+c^5+d^5-5*a*b*c*d*e+e^5));
     -- used 5.77 seconds

i3 : time V = integralClosure (R/I, Variable => X);
     -- used 4.23 seconds

i4 : ring ideal V

o4 = QQ [X ,X ,a,b,c,d,e, Degrees =>
          0  1                 {{5},{5},{1},{1},{1},{1},{1}}]

o4 : PolynomialRing

i5 : toString ideal V

o5 = ideal(a^2*b*c^2+b^2*c*d^2+a^2*d^2*e+a*b^2*e^2+c^2*d*e^2,
a*b^3*c+b*c^3*d+a^3*b*e+c*d^3*e+a*d*e^3,
a^5+b^5+c^5+d^5-5*a*b*c*d*e+e^5,
a*b*c^4-b^4*c*d-X_0*e-a^2*b^2*d*e+a*c^2*d^2*e+b^2*c^2*e^2
       -b*d^2*e^3,
a*b^2*c^3+X_1*d+a*b*c*d^2*e-a^2*b*e^3-d*e^5,
a^3*b^2*c-b*c^2*d^3-X_1*e-b^5*e-d^5*e+2*a*b*c*d*e^2,
a^4*b*c+X_0*d-a*b^4*e-2*b^2*c^2*d*e+a^2*c*d*e^2+b*d^3*e^2,
X_1*c+b^5*c+a^2*b^3*e-a*b*c^2*d*e-a*d^3*e^2,
X_0*c-a^2*b^2*c*d-b^2*c^3*e-a^4*d*e+2*b*c*d^2*e^2+a*b*e^4,
X_1*b-b*c^5+2*a*b^2*c*d*e-c^3*d^2*e+a^3*d*e^2-b*e^5,
X_0*b+a*b*c^2*d^2-b^3*c^2*e+a*d^4*e-a^2*b*c*e^2+b^2*d^2*e^2
     -c*d*e^4,
X_1*a-b^3*c^2*d+c*d^2*e^3,
X_0*a-b*c*d^4+c^4*d*e,
X_1^2+b^5*c^5+b^4*c^3*d^2*e+b*c^2*d^3*e^4+b^5*e^5+d^5*e^5,
X_0*X_1+b^3*c^4*d^3-b^2*c^7*e+b^2*c^2*d^5*e-b*c^5*d^2*e^2
       -a*b^2*c*d^3*e^3+b^4*c*d*e^4+a^2*b^2*d*e^5
       -a*c^2*d^2*e^5-b^2*c^2*e^6+b*d^2*e^7,
X_0^2+b*c^3*d^6+2*b^5*c*d^3*e+c*d^8*e-b^4*c^4*e^2
     +a^3*c^3*d^2*e^2+2*a^2*b^3*d^3*e^2-5*a*b*c^2*d^4*e^2
     +4*b^3*c^2*d^2*e^3-3*a*d^6*e^3+5*a^2*b*c*d^2*e^4
     -b^2*d^4*e^4-2*b*c^3*d*e^5-a^3*b*e^6+3*c*d^3*e^6-a*d*e^8)
\end{verbatim}

The command ``ring ideal V'' gives the polynomial ring and ``toString ideal V'' 
gives the defining ideal for the normalization.  

In {\sc Macaulay 2} we computed the radical of the three equations and then 
computed the normalization.  {\sc Singular} and {\sc Macaulay 2} use 
different algorithms for the command radical which means that there are 
examples which run well in one program, but not the other.  The radical 
computation done above in {\sc Macaulay 2} does not finish in {\sc Singular}
on the Pentium III, 600 MHZ machine with 256 MB of RAM, but this is a feature 
of this example.  We could just as easily find an example that {\sc Singular} 
completes and {\sc Macaulay 2} does not.  In {\sc Singular} we ran 
the normalization computation on the radical obtained from {\sc Macaulay 2}. 
Also this time we do not do the command show(nor);, because the output is 
similar to that of the first example.  Finally, as we did above, we alter the 
output a little to make it more readable for this paper.

\begin{verbatim}
> ring R = 0,(a,b,c,d,e),dp;
> ideal I = a2bc2+b2cd2+a2d2e+ab2e2+c2de2,
ab3c+bc3d+a3be+cd3e+ade3,
a5+b5+c5+d5-5abcde+e5,
a3b2cd-bc2d4+ab2c3e-b5de-d6e+3abcd2e2-a2be4-de6,
abc5-b4c2d-2a2b2cde+ac3d2e-a4de2+bcd2e3+abe5,
ab2c4-b5cd-a2b3de+2abc2d2e+ad4e2-a2bce3-cde5,
b6c+bc6+a2b4e-3ab2c2de+c4d2e-a3cde2-abd3e2+bce5,
a4b2c-abc2d3-ab5e-b3c2de-ad5e+2a2bcde2+cd2e4;
list nor = normal(I);

// 'normal' created a list of 1 ring(s).
// nor[1+1] is the delta-invariant in case of choose=wd.
// To see the rings, type (if the name of your list is nor):
     show( nor);
// To access the 1-st ring and map (similar for the others), 
     type: def R = nor[1]; setring R;  norid; normap;
// R/norid is the 1-st ring of the normalization and
// normap the map from the original basering to R/norid
> timer=1;
//used time: 3.21 sec
> def S = nor[1]; setring S; norid;
norid[1]=T(1)^2*T(2)*T(3)^2+T(2)^2*T(3)*T(4)^2
        +T(1)^2*T(4)^2*T(5)+T(1)*T(2)^2*T(5)^2
        +T(3)^2*T(4)*T(5)^2
norid[2]=T(1)*T(2)^3*T(3)+T(2)*T(3)^3*T(4)+T(1)^3*T(2)*T(5)
        +T(3)*T(4)^3*T(5)+T(1)*T(4)*T(5)^3
norid[3]=T(1)^5+T(2)^5+T(3)^5+T(4)^5
        -5*T(1)*T(2)*T(3)*T(4)*T(5)+T(5)^5
norid[4]=T(2)*T(3)*T(4)^4-T(3)^4*T(4)*T(5)-T(1)*T(6)
norid[5]=T(1)*T(2)*T(3)^2*T(4)^2-T(2)^3*T(3)^2*T(5)
        +T(1)*T(4)^4*T(5)-T(1)^2*T(2)*T(3)*T(5)^2
        +T(2)^2*T(4)^2*T(5)^2-T(3)*T(4)*T(5)^4+T(2)*T(6)
norid[6]=T(2)^3*T(3)^2*T(4)-T(3)*T(4)^2*T(5)^3-T(1)*T(7)
norid[7]=T(1)^2*T(2)^2*T(3)*T(4)+T(2)^2*T(3)^3*T(5)
+T(1)^4*T(4)*T(5)
        -2*T(2)*T(3)*T(4)^2*T(5)^2-T(1)*T(2)*T(5)^4-T(3)*T(6)
norid[8]=T(2)*T(3)^5-2*T(1)*T(2)^2*T(3)*T(4)*T(5)
+T(3)^3*T(4)^2*T(5)
        -T(1)^3*T(4)*T(5)^2+T(2)*T(5)^5-T(2)*T(7)
norid[9]=T(1)*T(2)*T(3)^4-T(2)^4*T(3)*T(4)
        -T(1)^2*T(2)^2*T(4)*T(5)+T(1)*T(3)^2*T(4)^2*T(5)
        +T(2)^2*T(3)^2*T(5)^2-T(2)*T(4)^2*T(5)^3-T(5)*T(6)
norid[10]=T(1)*T(2)^2*T(3)^3+T(1)*T(2)*T(3)*T(4)^2*T(5)
         -T(1)^2*T(2)*T(5)^3-T(4)*T(5)^5+T(4)*T(7)
norid[11]=T(2)^5*T(3)+T(1)^2*T(2)^3*T(5)
         -T(1)*T(2)*T(3)^2*T(4)*T(5)
         -T(1)*T(4)^3*T(5)^2+T(3)*T(7)
norid[12]=T(1)^3*T(2)^2*T(3)-T(2)*T(3)^2*T(4)^3-T(2)^5*T(5)
         -T(4)^5*T(5)+2*T(1)*T(2)*T(3)*T(4)*T(5)^2-T(5)*T(7)
norid[13]=T(1)^4*T(2)*T(3)-T(1)*T(2)^4*T(5)
         -2*T(2)^2*T(3)^2*T(4)*T(5)+T(1)^2*T(3)*T(4)*T(5)^2
         +T(2)*T(4)^3*T(5)^2+T(4)*T(6)
norid[14]=T(2)^4*T(3)*T(4)*T(5)^4+2*T(1)^3*T(4)^3*T(5)^4
         +4*T(1)^2*T(2)^2*T(4)*T(5)^5
         -4*T(1)*T(3)^2*T(4)^2*T(5)^5+2*T(2)*T(4)^2*T(5)^7
         +T(1)*T(2)*T(3)*T(4)*T(5)*T(6)
         +T(1)*T(3)^2*T(4)^2*T(7)-T(2)^2*T(3)^2*T(5)*T(7)
         -4*T(1)^2*T(3)*T(5)^2*T(7)-T(2)*T(4)^2*T(5)^2*T(7)
         +T(6)*T(7)
norid[15]=2*T(1)^2*T(2)*T(4)^4*T(5)^3
         -T(1)*T(2)^3*T(4)^2*T(5)^4-T(2)*T(3)^2*T(4)^3*T(5)^4
         +3*T(4)^5*T(5)^5+2*T(1)*T(2)*T(3)*T(4)*T(5)^6
         +3*T(2)^2*T(3)*T(4)^2*T(6)+2*T(1)*T(2)^2*T(5)^2*T(6)
         +2*T(3)^2*T(4)*T(5)^2*T(6)-2*T(3)^5*T(7)
         -3*T(4)^5*T(7)+3*T(1)*T(2)*T(3)*T(4)*T(5)*T(7)
         -2*T(5)^5*T(7)+2*T(7)^2
norid[16]=T(1)^3*T(2)*T(4)^5*T(5)-T(1)^3*T(3)^3*T(4)^2*T(5)^2
         -T(2)^2*T(4)^4*T(5)^4+T(2)*T(3)^3*T(4)*T(5)^5
         +T(1)^2*T(2)^2*T(4)*T(6)+T(2)^2*T(3)^2*T(5)*T(6)
         -2*T(2)*T(4)^2*T(5)^2*T(6)-T(6)^2
\end{verbatim}

{\sc Singular} used approximately 1 CPU second less than {\sc Macaulay 2} for 
this normalization computation.  As examples get larger this could 
be significant in obtaining success.
\end{example}

There are two other operations we want to illustrate.  Assume $A = R/I$ 
is an affine domain.  These first three examples give the integral 
closure of $A$ in the form $S/L$ where $S$ is of the form 
$R[Y_0,\ldots Y_t]$ and $S/L \cong \overline{A}$.  
If $A$ is a domain a user may want the fractions from the quotient 
field that generate the normalization of $A$ as a finite 
algebra over $A$.  This information can be obtained using both 
{\sc Macaulay 2} and {\sc Singular}.

\begin{example}\label{ex:fractions} 
For simplicity we use the same hypersurface as Example~\ref{ex1:hypersurface}.
\begin{verbatim}
i1 : A = QQ[x,y,z]/ideal(x^6-z^6-y^2*z^4);

i2 : time ICfractions(A)
     -- used 1.77 seconds

o2 = | x3/z2 x2/z x y z |

                  1            5
o2 : Matrix frac A  <--- frac A
\end{verbatim}

Thus $A[x^3/z^2,x^2/z] = \overline{A}$.  

In {\sc Singular} the method for obtaining the fractions is considerably 
different.  Several auxiliary functions are used.  The library normal.lib 
is no longer used and instead we need reesclos.lib.  {\sc Singular} 
is very particular how you define the ideal defining the 
affine domain.  It must be named ker or an error will be issued.  The output 
has the following form.
\begin{verbatim}
fraction;
[1]:
  xyz8
[2]:
  yz8
\end{verbatim}
Then the fraction we are looking for is [1] divided by [2], that is 
$\displaystyle{\frac{xyz^8}{yz^8} = x}$.  The fractions found in this way 
in {\sc Singular} are often not reduced.  

We use a new library in this example and since we included the 
printed output that is given when a library is loaded into {\sc Singular} 
in Example~\ref{ex1:hypersurface}, we removed it this time.
\begin{verbatim}
> LIB "reesclos.lib";
> ring R=0,(x,y,z),lp; ideal ker=x6-z6-y2z4;
> list L=primeClosure(R); closureRingtower(L);
> setring R(6); poly f=T(1); closureFrac(L); 
> setring R(1); fraction;
[1]:
   xyz8
[2]:
   yz8
> setring R(6); poly f=T(2); closureFrac(L); 
> setring R(1); fraction;
// ** redefining f **
[1]:
   y2z8
[2]:
   yz8
> setring R(6); poly f=T(3); closureFrac(L); 
> setring R(1); fraction;
// ** redefining f **
[1]:
   yz9
[2]:
   yz8
> setring R(6); poly f=T(4); closureFrac(L); 
> setring R(1); fraction;
// ** redefining f **
[1]:
   x2yz7
[2]:
   yz8
> setring R(6); poly f=T(5); closureFrac(L); 
> setring R(1); fraction;
// ** redefining f **
[1]:
   x3yz6
[2]:
   yz8
> timer=1;
//used time: 1.66 sec
\end{verbatim}
We can use this output from {\sc Macaulay 2} and {\sc Singular} to find a 
map from the variables used by Macaulay 2 to those used by Singular and 
vice versa.  
Also, given any element of the presentation of $\overline{A}$ as $S/J$, 
{\sc Singular} can find its image in the ring of fractions of $A$. This is 
not currently set up in {\sc Macaulay 2}.
\end{example}

Before computing the integral closure of an ideal, we 
point out one other subtle difference between the two programs.
In {\sc Macaulay 2} it is possible to define multigraded rings and 
the integralClosure program is designed to handle such rings, while 
in {\sc Singular} it is not possible to define such rings.  Multigraded 
rings play an 
important role in computing the integral closure of 
an ideal in {\sc Macaulay 2}.   
Since it is not possible to define a multigraded ring in {\sc Singular} the 
implementation there must be fundamentally different from that in 
{\sc Macaulay 2}.

It is possible, theoretically, to compute the integral closure 
of an ideal using both systems.  We say theoretically because the 
computation is often much too complex to complete, either due to 
memory or time (mostly memory).  The approach, in both systems, 
is the classical one, that is compute the Rees algebra first, 
find the normalization of the Rees algebra and then find the 
degree one piece of that graded algebra.  

In {\sc Macaulay 2} a function to compute $\overline{I}$ is not yet 
implemented in the main distribution, but we include one in the appendix 
which is being submitted to {\sc Macaulay 2} for inclusion.  
To run this program, a program to compute the Rees algebra that preserves 
the natural multigrading is needed.  We also include the code 
for this in the appendix.  There are several such programs for 
{\sc Macaulay 2} in circulation.  Besides the one in the appendix, one 
can be found in~\cite{M2book} and we 
have received yet another via personal communication from 
David~Eisenbud.  A package of such programs is being put together for 
inclusion with {\sc Macaulay 2}.  The one we include here is the only one that 
incorporates the natural multigrading, but is otherwise 
fundamentally the same as that in~\cite{M2book}, which could 
be easily altered to use the grading.  The program 
communicated by David~Eisenbud is more general and is based 
on his paper with Huneke and Ulrich~\cite{EHU}.  
Both computing blowups 
and integral closure of ideals are functions in {\sc Singular} if the Rees 
library has been loaded.  

We include two complete examples and discuss briefly two others.  The first 
example is 
one that can be checked by 
hand using standard techniques.  The second example is considerably 
more complicated and the remaining ideals illustrate how easily this process 
fails to complete.  All the examples are small in the sense that they use 
at most three generators, use at most three variables and have degree at most 
seven.

\begin{example}\label{simple}
Let $R = \QQ[x,y]$ and $I = (x^3,x^2y^2,y^7)$.  Using the 
fact that the integral closure of a monomial ideal corresponds 
to the integer lattice points in the convex hull of the exponent 
vectors of the ideal, this example is easily computed by hand.  
The first two commands load files named IdealNormal.m2 and 
blowup.m2 which are the programs given in the appendix.

\begin{verbatim}
i1 : load "IdealNormal.m2"; load "blowup.m2";
--loaded IdealNormal.m2
--loaded blowup.m2

i2 : R = QQ[x,y]; I = ideal(x^2,x*y^4,y^5);

o3 : Ideal of R
 
i4 : time idealIC(I)
     -- used 0.98 seconds

             2     3   5
o4 = ideal (x , x*y , y )

o4 : Ideal of R
\end{verbatim}

Computing the integral closure of an ideal in {\sc Singular} uses the 
library reesclos.lib which we used in Example~\ref{ex:fractions}.  The 
library must be loaded if it has not already been.
\begin{verbatim}
> ring R = 0,(x,y),dp; ideal I = x2,xy4,y5;
>  list J = normalI(I);  J;
[1]:
   _[1]=x2
   _[2]=y5
   _[3]=-xy3
\end{verbatim}

The output following [1]: is the generators of $\overline{I}$.

\end{example}

\begin{example}\label{complicated}
Let $R = \QQ[x,y,z]$ and $I = (y^6+x^2z,-x^6+y^4z^2)$. 
\begin{verbatim}
i5 : R = QQ[x,y,z, MonomialSize => 16]; 

i6 : I = ideal(y^6+x^2*z,-x^6+y^4*z^2); time idealIC(I)

o6 : Ideal of R 
     -- used 3.6 seconds

             6    2    6    4 2   5 4      2 3
o7 = ideal (y  + x z, x  - y z , x y  + x*y z )

o7 : Ideal of R
\end{verbatim}

\begin{verbatim}
> ring R=0,(x,y,z),dp; ideal I=y6+x2z,-x6+y4z2; 
> list J=normalI(I);  J;
[1]:
   _[1]=y6+x2z
   _[2]=-x6+y4z2
   _[3]=x5y4+xy2z3
\end{verbatim}
\end{example}

These examples show that both {\sc Macaulay 2} and {\sc Singular} can compute the 
integral closure of a 2 generated ideal in 3 variables, but an ideal with 
three generators will often cause problems.  
For example, if $R = \QQ[x,y,z]$ and $I = (y^6+x^2z,-x^6+y^4z^2,x^3y-z^3y)$ 
then {\sc Macaulay 2} will give a monomial overflow error and {\sc Singular} 
runs for about 7.1 hours and uses up the memory available on this machine.  
But, with two variables and two generators we can find a troublesome ideal 
in degree as low as seven.  If
$R = \QQ[x,y]$ and $I = (y^7+x^4,-y^6+xy^4)$ then {\sc Macaulay 2} 
finishes in 905.29 CPU seconds and uses about 66\% of the memory.  The 
ideal we get is 
$$\overline{I} = (x^2y^3+x^4,x^3y^2-x^4,x^4y+39x^3y^2-38x^4,x^5-x^4,y^6-xy^4,xy^5+x^4).$$
We ran this 
same example in {\sc Singular} and after 80,836.49 CPU seconds the machine 
quit the computation after using all of the memory.  Both programs 
struggle with various computations of ideal integral closure due to computing 
the integral closure of the blow-up ring which has more variables than 
the original ring.  Thus, while it is possible to compute the integral 
closure of some ideals 
using this method, for most cases it is impractical and the myriad of papers
dealing with special cases must be consulted.  

\section{Appendix}

The program idealIC, below, computes the integral closure of 
an ideal in {\sc Macaulay 2}.  This 
is the program used to compute the closures in 
Examples~\ref{simple},~\ref{complicated}. 
Except for ICfractionsLong all of the functions used in 
idealIC are self-explanatory, are commented, or are clearly explained 
in the {\sc Macaulay 2} documentation.  ICfractionsLong 
is similar to the function ICfractions 
used in Example~\ref{ex:fractions} which computes the 
normalization of $A$ and then returns a minimal generating set for 
$\overline{A}$ as an algebra over $A$. During the normalization computation 
fractions are computed that may be extraneous and while ICfractions only 
returns a minimal set, ICfractionsLong keeps track of all the fractions 
generated by the computation. To ensure that the entire degree one component of 
the integral closure of the Rees algebra is computed, ICfractionsLong is 
used rather than ICfractions.

\begin{verbatim}
idealIC = method()
idealIC(Ideal) := Ideal => (I) -> (
     R := ring I;
     n1 := numgens R;
     J1 := blowup(I);  --defining ideal of blow-up
     R2 := ring J1;
     n2 := numgens R2;
     S := R2/J1;  -- the blow-up
     Sfrac1 := first entries ICfractionsLong S;  
     -- The slow part is ICfractionsLong(S), as expected.
     Sfrac2 := apply(#Sfrac1-n2, i-> Sfrac1#i);
     toLift := select(Sfrac2, i-> (degree i)#1 == 1);
     -- toLift is the set of elements in the normalization  
     -- of the Rees algebra that generate the degree one 
     -- component.
     VarImage := flatten append (gens R,first entries gens I);
     -- We use VarImage to define the map from R2 -> R needed 
     -- to lift the fractions back to R.
     LiftMap := map(R,R2,VarImage);
     NewNums := apply(toLift, i-> LiftMap(
               substitute(numerator i, R2)));
     NewDenoms := apply(toLift, i-> LiftMap(
               substitute(denominator i, R2)));
     NewGens := apply(#toLift, i-> substitute(
               (NewNums#i)/(NewDenoms#i),R));
     ideal mingens (I + ideal(NewGens))
     )
\end{verbatim}

A program for computing the Rees algebra $R[tI]$ in {\sc Macaulay 2} which is used 
in the program IdealIC give above.

\begin{verbatim}
blowup = method(Options => {VarName => Y})
blowup(Ideal) := o-> (J)-> (
     -- Input:  J is any ideal.
     -- Output:  The result is the defining ideal of the 
     --          blowup algebra R[tJ] where R is the ring 
     --          of J. 
     -- METHOD:  We construct a polynomial ring 
     --          R[t,y_1,..,y_n] where n is the number of 
     --          generators of J.  We then construct the ideal 
     --          (Y_1-tJ_1,...,Y_n-tJ_n). Then R[tJ] is   
     --          isomorphic to the polynomial ring 
     --          R[Y_1,...,Y_n] mod the out put of this 
     --          algorithm.
     R := ring(J);
     n1 := numgens(J);
     n2 := numgens(R);
     Degs1 := flatten degrees source gens J; 
     -- Degs1 are the degrees of the generators of J.
     Degs1a := apply(n1,i->1+Degs1_i);
     Degs2 := flatten ((monoid R).Options.Degrees);
     -- Degs2 is the degrees of the variables of R.
     Degs := join({1},Degs2,Degs1a);
     S := coefficientRing(R)[t,gens R,
          (o.VarName)_(1)..(o.VarName)_(n1),
          MonomialOrder => Eliminate 1, Degrees => Degs];
     J1 := substitute(J,S);  -- Put J in the bigger ring S. 
     M := matrix{apply(n1,i-> S_(i+n2+1)-S_0*J1_i)};
     -- This is the matrix (Y_1-tJ_1,...,Y_n-tJ_n)
     K := gb(M); 
     -- Next eliminate t.
     L := ideal(selectInSubring(1,gens(K)));
     -- Set the Multidegrees.
     Degs1b := apply(n1,i->join({Degs1_i},{1}));
     Degs2b := apply(n2,i->join({Degs2_i},{0}));
     Degsb := join(Degs2b,Degs1b);
     Sb := coefficientRing(R)[gens R, 
           (o.VarName)_1..(o.VarName)_(n1),
           MonomialOrder => ProductOrder{n2,n1}, 
           Degrees => Degsb];
     trim(substitute(L,Sb))  
     )    
\end{verbatim}

\section*{Acknowledgements}
I thank Anna Guerrieri for encouraging me to write this 
article and Mike Stillman and Wolmer Vasconcelos for their helpful conversations. I 
thank David Eisenbud and C. Musili for giving me the opportunity to attend the 
conference in Hyderabad, India.


\end{document}